\documentclass[12pt]{article}
\usepackage[a4paper]{geometry}
\usepackage[latin9]{inputenc}
\usepackage{color}
\usepackage{calc}
\usepackage{amsmath}
\usepackage{amssymb}
\usepackage{graphicx}
\newtheorem{proposition}{Proposition}
\newtheorem{definition}{Definition}
\newtheorem{conjecture}{Conjecture}
\newtheorem{remark}{Remark}



\begin{document}

 \begin{center}{\Large \bf Elegant vertex labelings with prime numbers}\end{center}
 \vspace{0.2cm}
 
  \begin{center}{{\large Thierry Gensane}\\\medskip
  			 {\small LMPA J. Liouville\\
  			Universit\'e du Littoral\\
  			 Calais, FRANCE\\
  			\tt gensane@univ-littoral.fr\\}}\end{center}
  		 \vspace{0.2cm}
\begin{abstract}
We consider graph labelings with an assignment of odd prime numbers
to the vertices. Similarly to graceful graphs, a labeling is said
to be elegant if the absolute differences between the labels of adjacent
vertices describe exactly the first even numbers. The labels of an
elegant tree with $n$ vertices are the first $n$ odd prime numbers
and we want that the resulting edge labels are exactly the first even
numbers up to $2n-2$. We conjecture that each path is elegant and
we give the algorithm with which we got elegant paths of $n$ primes
for all $n$ up to $n=3500$. 
\end{abstract}
 		 \vspace{0.2cm}
\section{Introduction}

In this paper, we adapt the notion of graceful graphs by considering
an assignment of odd prime numbers to the vertices. Let $G=(V,E)$ be a graph,
we look for labelings of the vertices with distinct odd primes which induce
edge labelings with all even integers from $2$ up to $2|E|$. For
instance in the tree displayed in Fig.~\ref{fig_tree_12}, the first twelve odd
primes are assigned to the vertices and we get all the even positive integers up to
$22$. We call elegant any graph for which there exists such a labeling.
We refer to Galian~\cite{gali} for a very detailed survey about graph
labelings.

\begin{figure}[!t]
	\centering % Use \includegraphics to import figures; for example 
	\includegraphics[width=8.2cm,height=8.2cm]{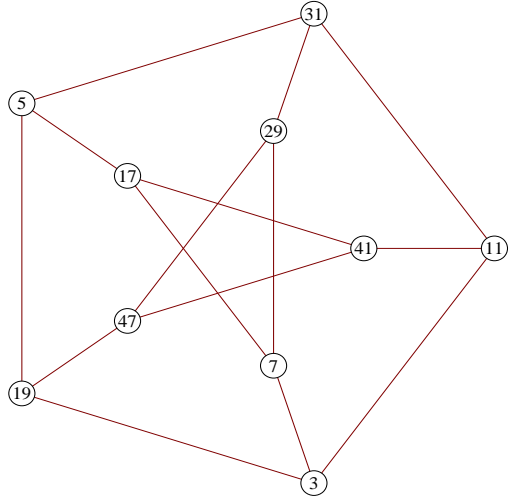}
	\vspace{-1cm}
	\caption{ \label{fig_petersen5}The Petersen graph is elegant.}
\end{figure}

Let us denote the increasing sequence
of all odd prime numbers by $p_{1},p_{2},\ldots,p_r,\ldots$ and $\mathbb{P}_{n}=\{p_{1},p_{2},\ldots,p_{n}\}$.
We now precise the definition of an elegant graph. 
\begin{definition} Let $G=(V,E)$ be an undirected graph without
loop or multiple edge, with $n$ vertices and $r$ edges. We say that
$G$ is {\rm elegant} if there exists an injective map $\varphi:V\to\left\{ p_{1},p_{2},\ldots,p_{r+1}\right\} $
such that the induced map 
\[
\begin{array}{rccl}
\psi\;: & E & \longrightarrow & 2\mathbb{N}^{*}\\
\\
 & e=uv & \longrightarrow & \psi(e)=|\varphi(v)-\varphi(u)|
\end{array}
\]
is a one-to-one correspondence from $E$ to $\{2,4,6,\ldots,2r\}$.

\end{definition}
\begin{center}
	\begin{figure}[!h]
		\begin{centering}
			\includegraphics{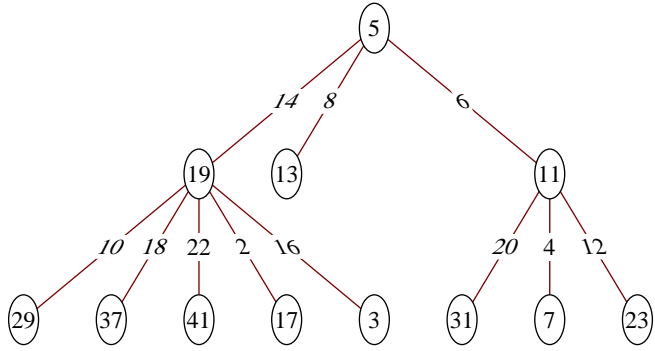} 
			\par\end{centering}
		\caption{\label{fig_tree_12} An elegant tree of $12$ vertices.}
	\end{figure}
	\par\end{center}
The complete graphs up to $K_{4}$  are
graceful and elegant: Their elegant labelings are respectively determined by $\varphi(V)=\{3,5\}, \{5,7,11\}$ and $\{7,11,17,19\}$.  As in the case of graceful graphs, it seems that no other elegant complete graphs exists. We display  an  elegant
labeling of the Petersen graph in Fig.~\ref{fig_petersen5}.

Of course, if the number of vertices is too weak relatively to the
maximal degree of the graph, then the graph is probably not (or cannot)
be elegant. For instance, the star graphs $S_{n}$ with $|V(S_{n})|=n$
are not elegant as soon as the center has more than $8$ adjacent
vertices (in fact, we verified that the only elegant star graphs are
$S_{2},$ $S_{3}$, $S_{5}$, $S_{6}$ and $S_{9}$). Nevertheless,
in the case of trees, we could hope that for each integer $d\geq2$,
the answer to the following question ${\rm A}_{d}$ be positive.

\bigskip{}

\noindent \textbf{Question ${\rm A}_{d}$:} Let $d\geq2$ be an integer.
Is there an integer $N_{d}$ such that each tree of maximal degree
$d$ and with more than $N_{d}$ vertices is elegant?

\bigskip{}

The answer to ${\rm A}_{d}$ is negative for $d$ large enough that is an easy consequence of the prime number theorem, see for instance
\cite{apos,hard}: Let us consider a symmetric and elegant tree $T$
rooted at a vertex of degree $d$, with exactly $d^{k}$ vertices
of degree $d+1$ at each level $k=1,\ldots,m-1$ and with $d^{m}$
leaves at level $m$. We have $n=|V(T)|=1+d+\cdots+d^{m}$ and then
$p_{n}\underset{d\to\infty}{\sim}md^{m}\log d$. Moreover, the absolute
difference between labels of two adjacent vertices is less than $2n-2\sim_{d\to\infty}2d^{m}$.
We consider the path $v_{0}v_{1}\cdots v_{k}\cdots v_{p}$ on the
tree $T$, from the vertex $v_{0}$ labeled by $\varphi(v_{0})=3$
and ended by the vertex $v_{p}$ labeled by $\varphi(v_{p})=p_{n}$.
Since $p\leq2m$, we get

\[
md^{m}\log d\underset{d\to\infty}{\sim}|\varphi(v_{0})-\varphi(v_{p})|\leq\sum_{i=0}^{^{p-1}}|\varphi(v_{i})-\varphi(v_{i+1})|\leq4md^{m}\left(1+o(1)\right)
\]

$\;$
%\vspace{-1.5cm}

%When $m=1$, we are
%back to the star graphs which cannot be elegant for %large $d$ as
%mentioned before.

%\medskip{}

%\begin{figure}[h]
%\centering % Use \includegraphics to import figures; for example 
% \includegraphics[width=5.2cm,height=5.2cm]{K3} \includegraphics[width=5.2cm,height=5.2cm]{K4}\includegraphics[width=5.2cm,height=5.2cm]{petersen5}
%\caption{ \label{fig_petersen5}The complete graphs $K_{3},$ $K_{4}$ and
%the Petersen graph are elegant.}
%\end{figure}

%\vspace{-1cm}
\noindent which is impossible when $d$ is large enough. Since $n\to\infty$
when $m\to\infty$, there exists a minimal integer $d_{1}$ such that
$A_{d}$ is false for all $d\geq d_{1}$.

As soon as $d=3$, it is difficult to know if $A_{d}$ is true or
false (and we are far from the value $d=e^{4}$ which appears in the previous proof). Let us illustrate this point
with an example: Let ${\rm C}_{n}$ be the regular caterpillar with
$n$ vertices of degree $3$ and $n+2$ leaves. On the one hand, up to
$n=25$ a trivial stocchastic algorithm has given elegant labelings
of $C_{n}$ only for $n=3,5,10,18,19,20,22$ but it is quite possible
that $C_{n}$ be elegant for all $n$ large enough. Let us recall
that Rosa \cite{rosa} proved that all caterpillars are graceful.
On the other hand, up to $n=25$,  as soon as we add a supplementary leaf $w$  on any leaf $v\in V(C_{n})$ 
 or if we supress
anywhere one leaf of $V(C_n)$, then surprisingly our program finds in a few seconds
that the modified tree is elegant.

%\begin{figure}[!h]
%\includegraphics{star9}
%\caption{\label{fig_star9}La plus grande étoile élégante.}
%\end{figure}

\medskip{}
 Fortunately, the case $d=2$ seems not to be resistant and we are
confident that $N_{2}=2$:

\begin{conjecture}For all $n\geq2$, the path of length $n$ is elegant.

\end{conjecture} Let us recall that a path of $n$ vertices is elegant
if there exists a permutation $\sigma\in S_{n}$ %on $\left\{ p_{1},p_{2},\ldots,p_{n}\right\} $
such that 
\[
\left\{ \left|p_{\sigma(i+1)}-p_{\sigma(i)}\right|;1\leq i\leq n-1\right\} =\left\{ 2,4,6,\ldots,2n-2\right\} .
\]

\noindent For instance, up to $n=10$, the following labelings are
elegant: 
\begin{itemize}
\item $3$ $\underset{2}{\text{\textendash\textendash}}$ $5$ 
\item $5$ $\underset{2}{\text{\textendash\textendash}}$ $3$ $\underset{4}{\text{\textendash\textendash}}$
$7$ 
\item $11$ $\underset{6}{\text{\textendash\textendash}}$ $5$ $\underset{2}{\text{\textendash\textendash}}$
3 $\underset{4}{\text{\textendash\textendash}}$ $7$ 
\item $13$ $\underset{6}{\text{\textendash\textendash}}$ $7$ $\underset{4}{\text{\textendash\textendash}}$
$11$ $\underset{8}{\text{\textendash\textendash}}$ $3$ $\underset{2}{\text{\textendash\textendash}}$
$5$ 
\item $7$ $\underset{2}{\text{\textendash\textendash}}$ $5$ $\underset{6}{\text{\textendash\textendash}}$
$11$ $\underset{8}{\text{\textendash\textendash}}$ $3$ $\underset{10}{\text{\textendash\textendash}}$
$13$ $\underset{4}{\text{\textendash\textendash}}$ $17$ 
\item  $13$ $\underset{4}{\text{\textendash\textendash}}$ $17$ $\underset{10}{\text{\textendash\textendash}}$
$7$ $\underset{12}{\text{\textendash\textendash}}$ $19$ $\underset{8}{\text{\textendash\textendash}}$
$11$ $\underset{6}{\text{\textendash\textendash}}$ $5$ $\underset{2}{\text{\textendash\textendash}}$
$3$ 
\item $13$ $\underset{10}{\text{\textendash\textendash}}$ $3$ $\underset{4}{\text{\textendash\textendash}}$
$7$ $\underset{2}{\text{\textendash\textendash}}$ $5$ $\underset{14}{\text{\textendash\textendash}}$
$19$ $\underset{8}{\text{\textendash\textendash}}$ $11$ $\underset{12}{\text{\textendash\textendash}}$
$23$ $\underset{6}{\text{\textendash\textendash}}$ $17$ 
\item $11$ $\underset{2}{\text{\textendash\textendash}}$ $13$ $\underset{8}{\text{\textendash\textendash}}$
$5$ $\underset{14}{\text{\textendash\textendash}}$ $19$ $\underset{10}{\text{\textendash\textendash}}$
$29$ $\underset{12}{\text{\textendash\textendash}}$ $17$ $\underset{6}{\text{\textendash\textendash}}$
$23$ $\underset{16}{\text{\textendash\textendash}}$ $7$ $\underset{4}{\text{\textendash\textendash}}3$ 
\item $19$ $\underset{12}{\text{\textendash\textendash}}$ $31$ $\underset{14}{\text{\textendash\textendash}}$
$17$ $\underset{6}{\text{\textendash\textendash}}$ $23$ $\underset{10}{\text{\textendash\textendash}}$
$13$ $\underset{16}{\text{\textendash\textendash}}$ $29$ $\underset{18}{\text{\textendash\textendash}}$
$11$ $\underset{8}{\text{\textendash\textendash}}$ $3$ $\underset{4}{\text{\textendash\textendash}}$
$7$ $\underset{2}{\text{\textendash\textendash}}$ $5$ 
\end{itemize}

\section{Operations on admissible paths}

In the sequel, a \textit{path of $l$ primes} ${\cal Q={\cal Q}}_{l}=q_{1}q_{2}\ldots q_{l}$
represents the labeling of the path of length $l$ with the primes
$q_{1},q_{2},\ldots,q_{l}$ in this order; we will  identify the vertices of a path and their label $q_i$.

\begin{definition} Let $n\geq2$ be a given integer. 
\begin{enumerate}
\item We say that a path ${\cal Q}_{l}=q_{1}q_{2}\ldots q_{l}$ is {\rm admissible}
if the primes $q_{i}\in\mathbb{P}_{n}$ are distinct and if the set
$E_{l}$ of the $l-1$ gaps $\vert q_{i+1}-q_{i}\vert$ is a subset
of cardinal $l-1$ of $\{2,4,\ldots,2n-2\}$. 
\item With regard to a path ${\cal Q}_{l}=q_{1}q_{2}\ldots q_{l}$ with $l<n$,
a prime $p$ in $\mathbb{P}_{n}$ is said to be {\rm free} if it
is not a vertex of ${\cal Q}_{l}$. A gap $2k\leq2n-2$ is said to
be {\rm free} if it does not belong to $E_{l}$. 
\end{enumerate}
\end{definition}

We also adopt the notations : 
\begin{itemize}
\item If ${\cal Q}=q_{1}q_{2}\ldots q_{l}$ then $\overline{{\cal Q}}=q_{l}q_{l-1}\ldots q_{1}$. 
\item If ${\cal Q}_{1}=q_{1}q_{2}\ldots q_{l}$ and ${\cal Q}_{2}=q_{l+1}q_{l+2}\ldots q_{p}$
, then ${\cal Q}_{1}{\cal Q}_{2}=q_{1}\ldots q_{l}q_{l+1}\ldots q_{p}$ is the     \textit{concatenation}  of ${\cal Q}_{1}$
and ${\cal Q}_{2}$. 
When we want to add a prime to a path ${\cal Q}={\cal Q}_{1}{\cal Q}_{2}$
respectively on the left end, between ${\cal Q}_{1}$ and ${\cal Q}_{2}$
or on the right end, we note these paths $p{\cal Q},$ ${\cal Q}_{1}p{\cal Q}_{2}$
or ${\cal Q}p$.
\item If ${\cal Q}=q_{1}q_{2}\ldots q_{l}$ then $f({\cal Q})=q_{1}$, $\ell({\cal Q})=q_{l}$
and ${\rm length}({\cal Q})=l$.
\item If ${\cal Q}=q_{1}q_{2}\ldots q_{l}$ is admissible,  $F_{p}({\cal Q})$ is
 \textit{the set of free primes for ${\cal Q}$} and $F_{g}({\cal Q})$ is  \textit{the set
of free gaps for ${\cal Q}$}. 
 
\end{itemize}

In the algorithm  described in Section 3, we  randomly apply transformations on admissible paths of length $l<n$ primes in order to find other admissible paths. Our aim is either to improve the length $l$ by adding a prime to the path, or to substitute a prime of the path for a free prime, or  simply to shuffle the path. Proposition \ref{shuffle} gives two elementary  tools $A1$-$A2$ and $A3$-$A4$ for shuffling an admissible path:

\begin{proposition}\label{shuffle} Let $n\geq 3$ and ${\cal Q}={\cal Q}_{1}{\cal Q}_{2}$
be an admissible path of $l<n$ primes. We denote by $\delta= |f({\cal Q}_{2})-\ell({\cal Q}_{1})|$ the gap between ${\cal Q}_{1}$ and ${\cal Q}_{2}$.
\begin{enumerate}\label{ 
	shuffle}
 
 \item[A1.] If  $|f({\cal Q}_{2})-f({\cal Q}_{1})|\in F_{g}({\cal Q})$, then the path ${\cal Q}^{*}=\overline{{\cal Q}}_{1}{\cal Q}_{2}$
is admissible and 
\[
F_{g}({\cal Q}^{*})=F_{g}({\cal Q})\cup\left\{ \delta\right\} \setminus\left\{ |f({\cal Q}_{2})-f({\cal Q}_{1})|\right\} .
\]
\item[A2.] If $|\ell({\cal Q}_{2})-\ell({\cal Q}_{1})|\in F_{g}({\cal Q})$, then the path ${\cal Q}^{*}={\cal Q}_{1}\overline{{\cal Q}}_{2}$
is admissible and 
\[
F_{g}({\cal Q}^{*})=F_{g}({\cal Q})\cup\left\{ \delta\right\} \setminus\left\{ |\ell({\cal Q}_{2})-\ell({\cal Q}_{1})|\right\} .
\]
\item[A3.] If $|\ell({\cal Q}_{2})-f({\cal Q}_{1})|\in F_{g}({\cal Q})$, then the path ${\cal Q}^{*}=Q_{2}Q_{1}$
is admissible and 
\[
F_{g}({\cal Q}^{*})=F_{g}({\cal Q})\cup\left\{ \delta\right\} \setminus\left\{ |\ell({\cal Q}_{2})-f({\cal Q}_{1})|\right\} .
\]
\item[A4.] If $|\ell({\cal Q}_{2})-f({\cal Q}_{1})|=\delta$, then the path ${\cal Q}^{*}=Q_{2}Q_{1}$ is admissible
and 
\[
F_{g}({\cal Q}^{*})=F_{g}({\cal Q}).
\]
\end{enumerate}
\end{proposition}

  \begin{remark}[Insertion of a free gap]\label{insertfreegap} {\rm In  Algorithm 1 decribded in Section 3 and when we  get an admissible transformation from  ${\cal Q}$ to  ${\cal Q}^*$ with a property $A_i$,  we  try  to insert a free prime in ${\cal Q}^*$: If $r$ is a new free prime for ${\cal  Q}^*$, we can try to insert it directly with Proposition~\ref{insertprime}. If a gap $ \delta$ has become free, we can try to insert in ${\cal Q}^*$ by Proposition~\ref{insertprime}, any prime $r=s\pm\delta\in F_p({\cal Q})$  if ${\cal Q}^*={\cal Q}_1{\cal Q}_2$ with  $s=\ell({\cal Q}_1)$ or $s=f({\cal Q}_2)$. Obviously, we also test in the algorithm if one of the paths $r{\cal Q}$ or ${\cal Q}r$ is admissible. 
  	
  	 For instance, let us consider $n=7$ and the admissible path  ${\cal Q}_6=7\textendash 19\textendash 17\textendash 11\textendash 3\textendash 13$. The last free gap is $4=|11-7|$ and we can apply the transformation given by A1 with  ${\cal Q}_1= 7\textendash 19\textendash 17$ and ${\cal Q}_2=11\textendash 3\textendash 13$, we get ${\cal Q}^{*}=\overline{{\cal Q}}_{1}{\cal Q}_{2}=17\textendash 19\textendash 7\textendash 11\textendash 3\textendash 13$ and the last free gap is now $\delta=6$. The last free prime is $r=5=11-\delta$ and we can insert $r$ by the third point of Proposition \ref{insertprime} with ${\cal Q}_1= 17\textendash 19\textendash 17\textendash 11$ and ${\cal Q}_2= 3\textendash 13$: We get an admissible and elegant path of seven primes with the admissible transformation $ {\cal Q}_6={\cal Q}_1{\cal Q}_2\to{\cal Q}_7={\cal Q}_{1}r\overline{{\cal Q}_{2}}=17\textendash 19\textendash 7\textendash 11\textendash 5 \textendash 13\textendash 3$.}
  	
  \end{remark}
  %\vspace{5cm}
%In order to insert a prime which becomes free after an admissible transformation of a %path of primes, we will use the following proposition.
\begin{proposition}[Insertion of a free prime] \label{insertprime}Let $n\geq3$
	and ${\cal Q=Q}_{1}{\cal Q}_{2}$ be an admissible path of $l<n$
	primes and let $r\in F_{p}({\cal Q})$.
We denote by \textbf{C}
	the condition 
	\[\textbf{C}\,:\,\left\{
	\begin{array}{c}
	\left(|r-p|\in F_{g}({\cal Q})\text{ and }|r-p|\neq|q-r|\in F_{g}({\cal Q})\right)\\
	\text{or}\\
	\left(|r-p|\in F_{g}({\cal Q})\text{ and }|q-r|=\delta\right)\\
	\text{or}\\
	\left(|q-r|\in F_{g}({\cal Q})\text{ and }|r-p|=\delta\right).
	\end{array}
	\right.
	\]
	
	 \noindent If  \textbf{C}
	is true for $\delta=|f({\cal Q}_2)-l({\cal Q}_1)|$ with
	%and $(p,q,r)$ defined by
	\begin{enumerate}
		\item $p=\ell({\cal Q}_{1})$ and $q=f({\cal Q}_{2})$, then the path ${\cal Q}^*=Q_{1}rQ_{2}$
		is admissible;
		\item $p=f({\cal Q}_{1})$ and $q=f({\cal Q}_{2})$, then the path ${\cal Q}^*=\overline{{\cal Q}_{1}}r{\cal Q}_{2}$
		is admissible;
		\item $p=\ell({\cal Q}_{1})$ and $q=\ell({\cal Q}_{2})$, then the path ${\cal Q}^*={\cal Q}_{1}r\overline{{\cal Q}_{2}}$
		is admissible.
	\end{enumerate}
\end{proposition}

\newpage 

In order to substitute a prime $p$ of a path ${\cal Q}$ for a free prime, we can consider $36$ transformations $A_5,A_6,\ldots,A_{40}:{\cal Q}\in X\rightarrow {\cal Q}^*\in Y$ 
 where $ X=\left\{ q{\cal Q}_{1}{\cal Q}_{2},{\cal Q}_{1}q{\cal Q}_{2},\right.$ $\left.{\cal Q}_{1}{\cal Q}_{2}q\right\} $
and
$
Y=$ $\left\{ r{\cal Q}_{1}{\cal Q}_{2},    r{\cal Q}_{1}{\cal \overline{Q}}_{2},     r{\cal \overline{Q}}_{1}{\cal Q}_{2},        r{\cal \overline{Q}}_{1}\overline{{\cal Q}}_{2},   {\cal Q}_{1}r{\cal Q}_{2},{\cal Q}_{1}r{\cal \overline{Q}}_{2},{\cal \overline{Q}}_{1}r{\cal Q}_{2},{\cal \overline{Q}}_{1}r\overline{{\cal Q}}_{2},\right.$ 

\noindent${\cal Q}_{1}{\cal Q}_{2}r,{\cal Q}_{1}{\cal \overline{Q}}_{2}r,
\left.{{\cal \overline{Q}}_{1}{\cal Q}_{2}r,\cal \overline{Q}}_{1}\overline{{\cal Q}}_{2}r\right\}.
$
 The following proposition
details only the transformations ${\cal Q}={\cal Q}_{1}q{\cal Q}_{2}$ $\to{\cal Q}^*={\cal Q}_{1}r{\cal Q}_2$
and ${\cal Q}={\cal Q}_{1}q{\cal Q}_{2}\to{\cal Q}^*={\cal \overline{Q}}_{1}r{\cal Q}_{2}$,
but it is trivial to find for which conditions the other transformations $A_i$ give admissible paths.

\begin{proposition} Let $n\geq 3$ and ${\cal Q}={\cal Q}_{1}q{\cal Q}_{2}$
be an admissible path of $l<n$ primes. Let $r$ be a free prime. 
\begin{enumerate}
%\item[A5.] 
\item The path ${\cal Q}^{*}={\cal Q}_{1}r{\cal Q}_{2}$ is admissible
if $ \left|r-\ell({\cal Q}_{1})\right|\neq \left|f({\cal Q}_{2})-r\right|$ and 
\[
\left\{ \left|r-\ell({\cal Q}_{1})\right|,\left|f({\cal Q}_{2})-r\right|\right\} \subset F_{g}(Q)\cup\left\{ \left|q-\ell({\cal Q}_{1})\right|,\left|f({\cal Q}_{2})-q\right|\right\} .
\]
%In that case, we have 
%\[
%F_{g}(Q^{*})=F_{g}(Q)\cup\left\{ \left|q-l({\cal Q}_{1})\right|,\left|f({\cal %Q}_{2})-q\right|\right\} \setminus\left\{ \left|p-l({\cal Q}_{1})\right|,\left|f({\cal %Q}_{2})-p\right|\right\} ,
%\]
%and $F_{p}(Q^{*})=F_{p}(Q)\cup\left\{ q\right\} \setminus\left\{ p\right\} $. 
\item The path ${\cal Q}^{*}={\cal \overline{Q}}_{1}r{\cal Q}_{2}$ is
admissible if $ \left|r-f({\cal Q}_{1})\right|\neq \left|f({\cal Q}_{2})-r\right|$ and
\[
\left\{ \left|r-f({\cal Q}_{1})\right|,\left|f({\cal Q}_{2})-r\right|\right\} \subset F_{g}(Q)\cup\left\{ \left|q-\ell({\cal Q}_{1})\right|,\left|f({\cal Q}_{2})-q\right|\right\} .
\]
%In that case, we have 
%\[
%F_{g}(Q^{*})=F_{g}(Q)\cup\left\{ \left|q-l({\cal Q}_{1})\right|,\left|f({\cal %Q}_{2})-q\right|\right\} \setminus\left\{ \left|p-l({\cal Q}_{1})\right|,\left|f({\cal %Q}_{2})-p\right|\right\} ,
%\]
%and $F_{p}(Q^{*})=F_{p}(Q)\cup\left\{ q\right\} \setminus\left\{ p\right\} $. 
\end{enumerate}
\end{proposition}

\section{An algorithm for finding elegant paths }

We have found elegant labelings for all paths up to $n=3500$ vertices
with the algorithm we  describe below. In Algorithm 1, we construct incrementally a sequence
of paths: from a path ${\cal Q}_{l}$ of $l$ primes, we try either
to add a prime or to modify the path ${\cal Q}_{l}$ without
improving its length. The integer $n$ being fixed, our aim is to
obtain $l=n$.

The  first part 1-4 of Algorithm 1 is a trivial greedy algorithm, we simply try to add free primes on the ends of ${\cal Q}_{l}$.
In the While statement of the step 5, we intensively and randomly use the transformations $A_i$ given in Section~2. When $m$ reaches $N$ without giving an elgant path of length $n$, we quit  and we  start a new run of Algorithm~1. With $N=40n$,
we get with Algorithm~1, all the elegant paths of $n$ primes from $n=2$ up to $n=200$
in less than 5 seconds with one core of a 3.6 GHz processor.

\medskip{}

If $n\geq200$, we accelerate the calculations with Algorithm 2 in
which we have chosen $N=20n$ et $c_{0}=20$: when a run reaches $l\geq n-\delta$,
we do not give up the path if $l<n$ but we suppress the right end $q_l$
of the path ${\cal Q}_{l}$. We take for instance $\delta=1$ for $n=200$ and $\delta=5$
from $n=2000$ up to $n=3500$. This can be seen as a perturbation
on a non-optimal and rigid configuration. When $n\in[2000,2100]$, the average run-time
to find one elegant path with Algorithm~2 is 4 minutes; when $n\in[3400,3500]$,
this average run-time becomes $40$ minutes. It is surprising to find  a solution so easily among $3500!>10^{10000}$ permutations.

\medskip{}

\newpage
\noindent %
\fbox{\begin{minipage}[t]{5.99 in}%
\textbf{Algorithm 1} %
\end{minipage}}

\noindent %
\fbox{\begin{minipage}[t]{5.99 in}%
		\begin{minipage}[t]{5.8in}
		\begin{enumerate}
			\item Randomly choose a prime $q_{1}$ in $\mathbb{P}_{n}$ and set $l:=1$,
			${\cal Q}_{l}:=q_{1}$; 
			\item $l_{1}:=0;$
			\item {\bf While} $l_{1}\neq l$ and $l<n$ {\bf do} 
			\begin{itemize}
			\item $l_{1}:=l_{1}+1$;
			
			\item If we find $p\in F_{p}({\cal Q}_{l})$ s.t. ${\cal Q}_{l}p$
			(or resp. $p{\cal Q}_{l}$) is admissible then  we set ${\cal Q}_{l+1}:={\cal Q}_{l}p$
			(or resp. ${\cal Q}_{l+1}:=p{\cal Q}_{l}$) and $l:=l+1$; 
			\end{itemize}
			\item {\bf If} $l=n$ {\bf then} return the elegant path ${\cal Q}_{n}$ and {\bf quit}.
		
		\end{enumerate}
	\end{minipage}
\end{minipage}}
\noindent %
\fbox{\begin{minipage}[t]{5.99 in}%
			\begin{minipage}[t]{5.8in}
\begin{enumerate}

\item[5.] {\bf While} $l<n$ and $m<N$ {\bf do} 
\begin{itemize}
	\item Randomly choose $i\in\{1,2,3\}$.
	
	{\bf Case} $i=1$ {\bf:} Look for an admissible transformation 
	${ \cal Q}^*=\overline{{ \cal Q}_{1}}{\cal Q}_{2}$  or ${ \cal Q}_{1}\overline{{\cal Q}_{2}}$ of  ${\cal Q}_{l}$ with A1 or A2.  If we succeed in changing ${\cal Q}_l$, we set ${\cal Q}_l:={\cal Q}^*$. Then,  we try  to insert a free prime in this new path ${\cal Q}_l$  with Remark~\ref{insertfreegap} and Prop.~\ref{insertprime}; if we succeed in this, we set $l:=l+1$ and  ${\cal Q}_l:=\cal Q^*$, the path  $\cal Q^*$ having been given by Prop. 6.
	
	{\bf Case} $i=2$ {\bf:} 	{\bf If}  $|\ell({\cal Q}_l)-f({\cal Q}_l)|\in F_g({\cal Q}_l)$ 	{\bf then} randomly choose $u\in\{1,\ldots,l-1\}$ and set ${\cal Q}_1:=q_1\cdots q_u$, ${\cal Q}_2:=q_{u+1}\cdots q_l$ and ${\cal Q}_l:={\cal Q}_2{\cal Q}_1$. Then,  if we succeed  in inserting a free prime in this new ${\cal Q}_l$  with Remark~\ref{insertfreegap} and Prop.~\ref{insertprime}, we set $l:=l+1$ and  ${\cal Q}_l:=\cal Q^*$, the path  $\cal Q^*$ having been given by Prop.~\ref{insertprime}.
	
		{\bf Else} there exists $u\in\{1,\ldots,l-1\}$ s.t. $|q_{u+1}-q_u|= |l({\cal Q}_l)-f({\cal Q}_l)|$
	and we set ${\cal Q}_l:={\cal Q}_2{\cal Q}_1$ with ${\cal Q}_1=q_1\cdots q_u$ and ${\cal Q}_2=q_{u+1}\cdots q_l$.
	
	{\bf Case} $i=3$ {\bf:} Look for a transformation among $A_5,\ldots,A_{40}$ which gives an admissible path ${\cal Q}^*$. If we succeed in modifying ${\cal Q}_l$, we set ${\cal Q}_l:={\cal Q}^*$.
	If we succeed  in inserting a free prime in this new ${\cal Q}_l$  with Remark~\ref{insertfreegap} and Prop.~\ref{insertprime}, we set $l:=l+1$ and  ${\cal Q}_l:=\cal Q^*$.

\item $m:=m+1$. 
\end{itemize}
\item[6.] {\bf Return} the path ${\cal Q}_{l}$ which is elegant if $l=n$.
\end{enumerate}
	\end{minipage}
 \end{minipage}}

\newpage

\noindent %
\fbox{\begin{minipage}[t]{5.99 in}%
\textbf{Algorithm 2} %
\end{minipage}}

\vspace{0.1cm}
 %Bermond and Schonheim [8] proved that all symmetrical
%trees are graceful.
\noindent %
\fbox{\begin{minipage}[t]{5.99 in}%
		\hspace{0.1in} \begin{minipage}[t]{5.8 in}
$l:=0$;

{\bf While} $l<n$ {\bf do}
\begin{enumerate}
\item {\bf Do} \textbf{Algorithm 1} which gives  an admissible
path ${\cal Q}_{l}$ of length $l$; 
\item {\bf If} $l=n$ {\bf then return} the elegant path ${\cal Q}_{n}$;
\item {\bf If} $n-\delta\leq l<n$ {\bf then}

\hspace{0.1cm} c:=0;

\hspace{0.1cm} {\bf While} $l<n$ and $c<c_{0}$ {\bf do} 
\begin{enumerate}
\item Suppress the last prime $q_{l}$ of ${\cal Q}_{l}$ and set $l:=l-1$; 
\item Do {\bf Step 5} of {\bf Algorithm 1} and set c:=c+1; 
\item {\bf If} $l=n$ {\bf then} return the elegant path ${\cal Q}_{n}$. 
\end{enumerate}
\end{enumerate}
	\end{minipage}
\end{minipage}}

\noindent \bigskip{}

Let us detail an example of a possible run of {\bf Algorithm~1} in
the case $n=11$. 
\begin{itemize}
\item {\bf Steps $1-4$} of  {\bf Algorithm 1}: We randomly choose $q_{1}=5$, $q_{2}=7,\;q_{3}=3,\ldots$, 
and we get 
\[
{\cal Q}_{10}:=5\underset{2}{\text{\textendash\textendash}}7\underset{4}{\text{\textendash\textendash}}3\underset{16}{\text{\textendash\textendash}}19\underset{12}{\text{\textendash\textendash}}31\underset{20}{\text{\textendash\textendash}}11\underset{18}{\text{\textendash\textendash}}29\underset{8}{\text{\textendash\textendash}}37\underset{14}{\text{\textendash\textendash}}23\underset{6}{\text{\textendash\textendash}}17,
\]
the last free prime is $13$ and the free gap is $10$. 
 \item {\bf Step $5$} of  {\bf Algorithm 1}:  We randomly choose  values for  $i$ in $\{1,2,3\}$  : 
\begin{itemize}
\item $i=1$: The gap $|17-7|=10$ is free and we  apply the transformation  ${\cal Q}={\cal Q}_{1}{\cal Q}_{2}\rightarrow {\cal Q}^*={\cal Q}_1\overline{{\cal Q}_{2}}$ with ${\cal Q}_{1}=5-7$ and  ${\cal Q}_{2}=3-19-31-11-29-37-23-17$. We get

\[
{\cal Q}_{10}:=
%{\cal Q}_1\overline{{\cal Q}_2}=
5\underset{2}{\text{\textendash\textendash}}7\underset{10}{\text{\textendash\textendash}}17\underset{6}{\text{\textendash\textendash}}23\underset{14}{\text{\textendash\textendash}}37\underset{8}{\text{\textendash\textendash}}29\underset{18}{\text{\textendash\textendash}}11\underset{20}{\text{\textendash\textendash}}31\underset{12}{\text{\textendash\textendash}}19\underset{16}{\text{\textendash\textendash}}3,
\]
the last free gap is now $4$. 
\item $i=3$ : We apply the transformation ${\cal Q}={\cal Q}_{1}{\cal Q}_{2}q\rightarrow {\cal Q}^*={\cal Q}_1r\overline{{\cal Q}_{2}}$ with ${\cal Q}_{1}=5-7-17$, ${\cal Q}_{2}=23-37-29-11-31-19$, $q=3$ and $r=13$. We find
\[
{\cal Q}_{10}:=5\underset{2}{\text{\textendash\textendash}}7\underset{10}{\text{\textendash\textendash}}17\underset{4}{\text{\textendash\textendash}}13\underset{6}{\text{\textendash\textendash}}19\underset{12}{\text{\textendash\textendash}}31\underset{20}{\text{\textendash\textendash}}11\underset{18}{\text{\textendash\textendash}}29\underset{8}{\text{\textendash\textendash}}37\underset{14}{\text{\textendash\textendash23}},
\]
the free prime is now $3$ and the free gap is $16$. 
\item$i=2$ : Since the gap between the two ends is $|23-5|=18\notin F_g({\cal Q})$,  we apply  ${\cal Q}={\cal Q}_{1}{\cal Q}_{2}\rightarrow {\cal Q}^*={\cal Q}_2{\cal Q}_{1}$ with ${\cal Q}_{1}=5-7-17-13-19-31-11$ and  ${\cal Q}_{2}=29-37- 23$. We get
\[
{\cal Q}_{10}:=29\underset{8}{\text{\textendash\textendash}}37\underset{14}{\text{\textendash\textendash}}23\underset{18}{\text{\textendash\textendash}}5\underset{2}{\text{\textendash\textendash}}7\underset{10}{\text{\textendash\textendash}}17\underset{4}{\text{\textendash\textendash}}13\underset{6}{\text{\textendash\textendash}}19\underset{12}{\text{\textendash\textendash}}31\underset{20}{\text{\textendash\textendash}}11.
\]
\item $i=1$ : The gap $|29-13|=16$ is free and we  apply the transformation  ${\cal Q}={\cal Q}_{1}{\cal Q}_{2}\rightarrow {\cal Q}^*=\overline{{\cal Q}_1}{\cal Q}_{2}$ with ${\cal Q}_{1}=29-37-23-5-7$ and  ${\cal Q}_{2}=17-13-19-31-11$. We get
\[
{\cal Q}_{10}:=17\underset{10}{\text{\textendash\textendash}}7\underset{2}{\text{\textendash\textendash5}}\underset{18}{\text{\textendash\textendash}}23\underset{14}{\text{\textendash\textendash}}37\underset{8}{\text{\textendash\textendash}}29\underset{16}{\text{\textendash\textendash}}13\underset{6}{\text{\textendash\textendash}}19\underset{12}{\text{\textendash\textendash}}31\underset{20}{\text{\textendash\textendash}}11,
\]
and the last free gap is $4$. We can place the prime $3$ between
$7$ and $5$ and we get an elegant path of eleven primes: 
\[
{\cal Q}_{11}:=17\underset{10}{\text{\textendash\textendash}}7\underset{4}{\text{\textendash\textendash3}}\underset{2}{\text{\textendash\textendash5}}\underset{18}{\text{\textendash\textendash}}23\underset{14}{\text{\textendash\textendash}}37\underset{8}{\text{\textendash\textendash}}29\underset{16}{\text{\textendash\textendash}}13\underset{6}{\text{\textendash\textendash}}19\underset{12}{\text{\textendash\textendash}}31\underset{20}{\text{\textendash\textendash}}11.
\]
\end{itemize}
\end{itemize}
%%%%%%%%%%%%%%%%%%%%%%%%%%%%%%%%%%%%%%%%%%%%%%%%%%%%%%%

\subsection*{Acknowledgements}

Thanks to  Shalom Eliahou for his encouragment and more. %Lemma~\ref{lem:Technical}.

%%%%%%%%%%%%%%%%%%%%%%%%%%%%%%%%%%%%%%%%%%%%%%%%%%%%%%%
% You do not have to use the same format for your references, but 
%    include everything in this file.  Don't use natbib please.
% If you use BibTeX to create a bibliography, copy the .bbl file into here.

\newpage
\end{document}